# New Wallis- and Catalan-Type Infinite Products for $\pi$, $e$, and $\sqrt{2+\sqrt{2}}$

## Jonathan Sondow and Huang Yi

**Abstract**

We generalize Wallis's 1655 infinite product for $\pi/2$ to one for $(\pi/K)\csc(\pi/K)$, as well as give new Wallis-type products for $\pi/4$, $2$, $\sqrt{2+\sqrt{2}}$, $2\pi/3\sqrt{3}$, and other constants. The proofs use a classical infinite product formula involving the gamma function. We also extend Catalan's 1873 infinite product of radicals for $e$ to Catalan-type products for $e/4$, $\sqrt{e}$, and $e^{3/2}/2$. Here the proofs use Stirling's formula. Finally, we find an analog for $e^{2/3}/\sqrt{3}$ of Pippenger's 1980 product for $e/2$, and conjecture that they can be generalized to a product for a power of $e^{1/K}$.

**1. INTRODUCTION.** In 1655 Wallis [16] published his famous infinite product for pi:

$$\frac{\pi}{2} = \frac{2}{1}\frac{2}{3}\frac{4}{3}\frac{4}{5}\frac{6}{5}\frac{6}{7}\frac{8}{7}\frac{8}{9}\frac{10}{9}\frac{10}{11}\frac{12}{11}\frac{12}{13}\frac{14}{13}\frac{14}{15}\frac{16}{15}\cdots. \quad (1)$$

In 1873 Catalan [4] proved the Wallis-type formulas

$$\frac{\pi}{2\sqrt{2}} = \frac{4}{3}\frac{4}{5}\frac{8}{7}\frac{8}{9}\frac{12}{11}\frac{12}{13}\frac{16}{15}\cdots \quad (2)$$

and

$$\sqrt{2} = \frac{2}{1}\frac{2}{3}\frac{6}{5}\frac{6}{7}\frac{10}{9}\frac{10}{11}\frac{14}{13}\frac{14}{15}\cdots. \quad (3)$$

Together they give a beautiful factorization of Wallis's formula, which we write symbolically as $(1) = (2) \times (3)$.

Catalan also found a product for *e* similar to Wallis's product for pi:

$$e = \frac{2}{1}\left(\frac{4}{3}\right)^{1/2}\left(\frac{6}{5}\frac{8}{7}\right)^{1/4}\left(\frac{10}{9}\frac{12}{11}\frac{14}{13}\frac{16}{15}\right)^{1/8}\cdots. \quad (4)$$

Even closer to Wallis's formula is a product for *e* discovered by Pippenger [12] in 1980:

$$\frac{e}{2} = \left(\frac{2}{1}\right)^{1/2} \left(\frac{2}{3}\frac{4}{3}\right)^{1/4} \left(\frac{4}{5}\frac{6}{5}\frac{6}{7}\frac{8}{7}\right)^{1/8} \left(\frac{8}{9}\frac{10}{9}\frac{10}{11}\frac{12}{11}\frac{12}{13}\frac{14}{13}\frac{14}{15}\frac{16}{15}\right)^{1/16} \cdots . \tag{5}$$

While Catalan obtained (2), (3), and (4) as by-``products'' of series and integrals for the gamma function $\Gamma(x)$, Pippenger's proof of (5) uses only Stirling's asymptotic formula

$$N! \sim \sqrt{2\pi N}(N/e)^N \quad (N \to \infty).$$

In this note, we offer several new products like Wallis's.

**Theorem 1.** *The following Wallis-type formulas are valid:*

$$\frac{\pi}{4} = \frac{2}{3}\frac{6}{5}\frac{8}{7}\frac{8}{9} \cdot \frac{10}{11}\frac{14}{13}\frac{16}{15}\frac{16}{17} \cdot \frac{18}{19}\frac{22}{21}\frac{24}{23}\frac{24}{25} \cdots , \tag{6}$$

$$2 = \frac{2}{1}\frac{4}{3}\frac{4}{5}\frac{6}{7} \cdot \frac{10}{9}\frac{12}{11}\frac{12}{13}\frac{14}{15} \cdot \frac{18}{17}\frac{20}{19}\frac{20}{21}\frac{22}{23} \cdots , \tag{7}$$

$$\frac{\pi}{4\sqrt{2-\sqrt{2}}} = \frac{8}{7}\frac{8}{9}\frac{16}{15}\frac{16}{17}\frac{24}{23}\frac{24}{25} \cdots , \tag{8}$$

$$\sqrt{2-\sqrt{2}} = \frac{2}{3}\frac{6}{5}\frac{10}{11}\frac{14}{13}\frac{18}{19}\frac{22}{21} \cdots , \tag{9}$$

$$\frac{\pi}{2\sqrt{2+\sqrt{2}}} = \frac{2}{3}\frac{4}{3}\frac{4}{5}\frac{6}{5}\frac{8}{7}\frac{8}{9} \cdot \frac{10}{11}\frac{12}{11}\frac{12}{13}\frac{14}{13}\frac{16}{15}\frac{16}{17} \cdot \frac{18}{19}\frac{20}{19}\frac{20}{21}\frac{22}{21}\frac{24}{23}\frac{24}{25} \cdots , \tag{10}$$

$$\sqrt{2+\sqrt{2}} = \frac{2}{1}\frac{6}{7}\frac{10}{9}\frac{14}{15}\frac{18}{17}\frac{22}{23} \cdots , \tag{11}$$

$$\frac{2\pi}{3\sqrt{3}} = \frac{3}{2}\frac{3}{4}\frac{6}{5}\frac{6}{7}\frac{9}{8}\frac{9}{10} \cdots . \tag{12}$$

*Moreover, Wallis's product and* (12) *are the cases $K=2$ and $3$ of the following general formula, which holds for any integer $K \geq 2$:*

$$\frac{\pi/K}{\sin(\pi/K)} = \frac{K}{K-1}\frac{K}{K+1}\frac{2K}{2K-1}\frac{2K}{2K+1}\frac{3K}{3K-1}\frac{3K}{3K+1} \cdots . \tag{13}$$

Theorem 1 evidently gives the factorizations of Wallis's formula

$$(1) = (6) \times (7) = (8) \times (9) \times (7) = (10) \times (11) \tag{14}$$

as well as the factorization (3) = (9) × (11) of Catalan's product for

$$\sqrt{2} = \sqrt{2 - \sqrt{2}}\sqrt{2 + \sqrt{2}}.$$

In particular, formulas (1), (3), (7), and (9) imply (6), (8), (11), and (10), by division.

On the other hand, the product (7) for 2 cannot be proved by simply squaring Catalan's product (3) for $\sqrt{2}$. Much less can the square of the product (11) for $\sqrt{2 + \sqrt{2}}$ be obtained from (7) and (3) by addition!

We also give some new products of the same shape as Catalan's for $e$.

**Theorem 2.** *The following Catalan-type formulas hold:*

$$\frac{e}{4} = \left(\frac{2}{3}\right)^{1/2}\left(\frac{4\;6}{5\;7}\right)^{1/4}\left(\frac{8\;10\;12\;14}{9\;11\;13\;15}\right)^{1/8}\cdots, \tag{15}$$

$$\sqrt{e} = \frac{2}{1}\left(\frac{2}{3}\right)^{1/2}\left(\frac{6\;6}{5\;7}\right)^{1/4}\left(\frac{10\;10\;14\;14}{9\;11\;13\;15}\right)^{1/8}\cdots, \tag{16}$$

$$\frac{e^{3/2}}{2} = \frac{2}{1}\left(\frac{4}{3}\right)^{1/2}\left(\frac{4\;8}{5\;7}\right)^{1/4}\left(\frac{8\;12\;12\;16}{9\;11\;13\;15}\right)^{1/8}\cdots, \tag{17}$$

$$\frac{e^{2/3}}{\sqrt{3}} = \left(\frac{3}{2}\right)^{1/3}\left(\frac{3\;6\;6\;9}{4\;5\;7\;8}\right)^{1/9}\left(\frac{9\;12\;12\;15\;15\;18\;18\;21\;21\;24\;24\;27}{10\;11\;13\;14\;16\;17\;19\;20\;22\;23\;25\;26}\right)^{1/27}\cdots. \tag{18}$$

Note that Pippenger's, Catalan's, and our first three products for $e$ are related by the factorizations

$$(5)^2 = (4) \times (15) = \frac{1}{2}(16) \times (17). \tag{19}$$

Also, notice that the ``geometric-series product''

$$2 = 2^{1/2}2^{1/4}2^{1/8}\cdots \tag{20}$$

allows the factorization (4) = (15) × (20)$^2$. Thus Catalan's and Pippenger's products are related by the surprisingly simple formula (4) = (5) × (20).

We built the product (18) for $e^{2/3}/\sqrt{3}$ from the product (12) for $2\pi/3\sqrt{3}$ by analogy with Pippenger's construction of his formula (5) for $e/2$ from Wallis's formula (1) for $\pi/2$. Just as (13) generalizes (1) and (12), so too *we conjecture that one can generalize* (5) *and* (18) *to a product for a power of $e^{1/K}$*.

Formulas (7), (3), and (1) for 2, $\sqrt{2}$, and $\pi$ show that a Wallis-type infinite product (in particular, a factor of Wallis's product) can be rational, algebraic irrational, or transcendental, respectively. As for Catalan-type infinite products, we know of no algebraic example; in fact, the known ones are all algebraically dependent on $e$.

Other types of infinite products can be found for $\pi$ in [**2**, p. 94], [**3**, p. 55], [**5**, Section 1.4.2], [**7**], [**11**], [**19**], for $e$ in [**6**], [**8**], [**18**], and for both $\pi$ and $e$ in [**9**], [**10**], [**13**].

We conclude the Introduction with some remarks on the use of modern symbolic algorithms, as implemented in computer algebra systems such as *Mathematica*. They can compute symbolically not only particular Wallis-type products such as (1) and (12), but also general ones like (13). This is useful for generating new Wallis-type formulas. On the other hand, *Mathematica* (version 7.0.0) cannot evaluate the Catalan-type infinite products symbolically. (Surprisingly, neither can it calculate them numerically. The reason given is ``overflow.'') In any case, such a computer-assisted approach is only intended for mechanical discovery and does not provide a rigorous mathematical proof. Thus the method is advisory, and provides us with a modern scope within which we can do more mathematics.

**2. PROOF OF THEOREM 1.** In view of the factorizations (14) and the fact that (12) is a special case of (13), it suffices to prove (7), (9), and (13). We use the following classical formula [**17**, Section 12.13], which is a corollary of the Weierstrass infinite product for the gamma function.

*If $k$ is a positive integer and $a_1 + a_2 + \cdots + a_k = b_1 + b_2 + \cdots + b_k$, where the $a_j$ and $b_j$ are complex numbers and no $b_j$ is zero or a negative integer, then*

$$\prod_{n=0}^{\infty} \frac{(n+a_1)\cdots(n+a_k)}{(n+b_1)\cdots(n+b_k)} = \frac{\Gamma(b_1)\cdots\Gamma(b_k)}{\Gamma(a_1)\cdots\Gamma(a_k)}. \tag{21}$$

To prove (7), we write the product as

$$\prod_{n=0}^{\infty} \frac{(8n+2)(8n+4)(8n+4)(8n+6)}{(8n+1)(8n+3)(8n+5)(8n+7)} = \prod_{n=0}^{\infty} \frac{(n+(1/4))(n+(1/2))(n+(1/2))(n+(3/4))}{(n+(1/8))(n+(3/8))(n+(5/8))(n+(7/8))}.$$

Then from (21) and Euler's reflection formula [**17**, Section 12.14]

$$\Gamma(x)\Gamma(1-x) = \frac{\pi}{\sin \pi x}, \tag{22}$$

the product is equal to

$$\frac{\Gamma(1/8)\Gamma(3/8)\Gamma(5/8)\Gamma(7/8)}{\Gamma(1/4)\Gamma(1/2)\Gamma(1/2)\Gamma(3/4)} = \frac{\sin(\pi/4)\sin(\pi/2)}{\sin(\pi/8)\sin(3\pi/8)} = \frac{1/\sqrt{2}}{\left(1/\sqrt{4+\sqrt{8}}\right)\left(1/\sqrt{4-\sqrt{8}}\right)} = 2,$$

as claimed. (For the values of sine, see for example [**15**, Chap. 32].)

Similarly, the product in (9) can be written

$$\frac{2 \cdot 6}{3 \cdot 5} \frac{10 \cdot 14}{11 \cdot 13} \frac{18 \cdot 22}{19 \cdot 21} \cdots = \prod_{n=0}^{\infty} \frac{(8n+2)(8n+6)}{(8n+3)(8n+5)} = \prod_{n=0}^{\infty} \frac{(n+(1/4))(n+(3/4))}{(n+(3/8))(n+(5/8))},$$

and by (21) and (22) its value is

$$\frac{\Gamma(3/8)\Gamma(5/8)}{\Gamma(1/4)\Gamma(3/4)} = \frac{\sin(\pi/4)}{\sin(3\pi/8)} = \frac{1/\sqrt{2}}{1/\sqrt{4-\sqrt{8}}} = \sqrt{2-\sqrt{2}},$$

as desired.

Finally, to prove (13) we write the product as

$$\prod_{n=1}^{\infty} \frac{(Kn)^2}{(Kn-1)(Kn+1)} = \prod_{n=0}^{\infty} \frac{(K(n+1))^2}{(K(n+1)-1)(K(n+1)+1)} = \prod_{n=0}^{\infty} \frac{(n+1)^2}{(n+1-(1/K))(n+1+(1/K))}.$$

Then by (21), (22), and the factorial property $\Gamma(x+1) = x\Gamma(x)$ (see [**17**, Section 12.12]), together with the value $\Gamma(1) = 0! = 1$, the product is equal to

$$\frac{\Gamma(1-(1/K))\Gamma(1+(1/K))}{\Gamma(1)^2} = \Gamma\left(1-\frac{1}{K}\right)\frac{1}{K}\Gamma\left(\frac{1}{K}\right) = \frac{\pi/K}{\sin(\pi/K)},$$

as was to be shown. This completes the proof of Theorem 1. •

In a similar way, one can also derive Catalan's products (2) and (3) for $\pi/2\sqrt{2}$ and $\sqrt{2}$ (compare [**1**, Section 4.4]). For a different way to prove (2), see [**14**, Example 7].

**3. PROOF OF THEOREM 2.** By the factorizations (19), to prove the first three formulas it suffices to prove (16). Equivalently, we show that

$$\left(\frac{2}{3}\right)^{1/2} \left(\frac{6 \cdot 6}{5 \cdot 7}\right)^{1/4} \left(\frac{10 \cdot 10 \cdot 14 \cdot 14}{9 \cdot 11 \cdot 13 \cdot 15}\right)^{1/8} \cdots = \frac{\sqrt{e}}{2}. \tag{23}$$

Notice the cancellations that occur in computing the partial products, which are simply

$$\frac{2^{1/2}}{3^{1/2}}, \frac{2^{2/2}}{(5 \cdot 7)^{1/4}}, \frac{2^{3/2}}{(9 \cdot 11 \cdot 13 \cdot 15)^{1/8}}, \cdots, \frac{2^{n/2}}{\left((2^n+1)(2^n+3)\cdots(2^{n+1}-1)\right)^{1/2^n}}, \cdots.$$

Using factorials, we can write the $n$th partial product as

$$2^n!2\left(\frac{2^n!(2^n+2)(2^n+4)\cdots 2^{n+1}}{2^{n+1}!}\right)^{1/2^n} = 2^n!2\left(\frac{2^n!2^{2^{n-1}}(2^{n-1}+1)(2^{n-1}+2)\cdots 2^n}{2^{n+1}!}\right)^{1/2^n}$$

$$= 2^{\frac{n+1}{2}}\left(\frac{(2^n!)^2}{2^{n-1}!2^{n+1}!}\right)^{1/2^n}.$$

After applying Stirling's formula to each factorial, the resulting expression simplifies to $\sqrt{e}/2$, establishing (23).

To evaluate the infinite product in (18), note first that, for $k \geq 2$, its $k$th factor is the $1/3^k$ power of the product of $2(3^{k-1} - 3^{k-2}) = 4 \cdot 3^{k-2}$ fractions. Factoring powers of 3 out of their numerators, we see that we can write the infinite product as

$$\left(\frac{3^1}{2}\right)^{1/3}\left(\frac{2^2 3^5}{4\cdot 5\cdot 7\cdot 8}\right)^{1/9}\left(\frac{(2\cdot 4\cdot 5\cdot 7\cdot 8)^2 3^{17}}{10\cdot 11\cdot 13\cdot 14\cdot 16\cdot 17\cdot 19\cdot 20\cdot 22\cdot 23\cdot 25\cdot 26}\right)^{1/27}\cdots,$$

and that by induction the exponents $1, 5, 17, \ldots$ are equal to

$$2\cdot 3^{k-1} - 1 = 4\cdot 3^{k-2} + 2\cdot 3^{k-2} - 1,$$

for $k = 1, 2, 3, \ldots$. Then the partial products are

$$\left(\frac{3^1}{2}\right)^{1/3}, \left(\frac{3^8}{2\cdot 4\cdot 5\cdot 7\cdot 8}\right)^{1/9}, \left(\frac{3^{41}}{2\cdot 4\cdot 5\cdot 7\cdot 8\cdot 10\cdot 11\cdot 13\cdot 14\cdot 16\cdot 17\cdot 19\cdot 20\cdot 22\cdot 23\cdot 25\cdot 26}\right)^{1/27}, \ldots$$

where the exponents $1, 8, 41, \ldots$ are given by the formula

$$E_n := \sum_{k=1}^{n} 3^{n-k}(2\cdot 3^{k-1} - 1) = 2n3^{n-1} - \frac{1}{2}(3^n - 1),$$

for $n = 1, 2, 3, \ldots$. Thus the $n$th partial product is

$$\left(\frac{3^{E_n}}{2\cdot 4\cdot 5\cdot 7\cdots(3^n - 4)(3^n - 2)(3^n - 1)}\right)^{1/3^n} = \left(\frac{3^{E_n}}{3^n!/(3^{3^{n-1}}3^{n-1}!)}\right)^{1/3^n} = 3^{1/3}\left(\frac{3^{E_n}3^{n-1}!}{3^n!}\right)^{1/3^n},$$

which by Stirling's formula is asymptotic to $e^{2/3}/3^{1/2}$. This proves (18) and completes the proof of Theorem 2. ∙

**ACKNOWLEDGMENT.** The authors are grateful to the anonymous referee who pointed out formulas (12) and (13) and whose helpful comments prompted our remarks on the use of computers.

*209 West 97th Street, New York, NY 10025*
*jsondow@alumni.princeton.edu*

*Department of Mathematics, Beijing Normal University, Beijing 100875, P. R. China*
*peiweihuang@mail.bnu.edu.cn*